\documentclass[11pt]{amsart}
\usepackage[T1]{fontenc}
\usepackage[utf8]{inputenc}
\usepackage{lmodern}
\usepackage{amsmath,amssymb,amsthm,mathtools}
\usepackage{microtype}
\usepackage[margin=1.12in]{geometry}
\usepackage[hidelinks]{hyperref}
\usepackage{enumitem}
\setlist[itemize]{leftmargin=*,topsep=3pt,itemsep=3pt}
\allowdisplaybreaks
\newtheorem{theorem}{Theorem}[section]
\newtheorem{proposition}[theorem]{Proposition}
\newtheorem{corollary}[theorem]{Corollary}
\newtheorem{lemma}[theorem]{Lemma}
\theoremstyle{definition}

\theoremstyle{remark}
\newtheorem{remark}[theorem]{Remark}
\newtheorem{example}[theorem]{Example}
\newcommand{\CP}{\mathbb{CP}}
\newcommand{\CC}{\mathbb C}
\newcommand{\RR}{\mathbb R}

\newcommand{\Ric}{\operatorname{Ric}}
\newcommand{\Kill}{\mathfrak{kill}}
\newcommand{\PU}{\operatorname{PU}}
\newcommand{\FS}{\mathrm{FS}}
\newcommand{\Lie}{\operatorname{Lie}}
\newcommand{\id}{\mathrm{id}}

\title[Positivity for projectively induced K\"ahler--Einstein metrics]{A local symmetry forces positivity for \mbox{projectively induced K\"ahler--Einstein metrics}}
\author{Filippo Salis}
\address{Dipartimento di Scienze Matematiche, Politecnico di Torino, Corso Duca degli Abruzzi 24, 10129 Torino, Italy}
\email{filippo.salis@polito.it}
\date{\today}
\begin{document}
\begin{abstract}
We prove that a projectively induced K\"ahler--Einstein metric has positive Einstein constant as soon as it admits one nonzero local real-holomorphic Killing field. Neither completeness nor a prescribed global group action is assumed. This gives a partial answer, under a local symmetry hypothesis, to the general positivity problem for K\"ahler--Einstein submanifolds of (finite-dimensional) projective space. The key observation is that, by Calabi rigidity, a local real-holomorphic Killing field is induced by a projective unitary vector field whose tangency propagates along any connected global K\"ahler immersion, without a simple-connectivity assumption. Applied to Hulin's complete continuation, this produces a global Killing field with a bounded Hamiltonian potential given by the ambient projective moment map. The resulting Laplace equation, the Omori--Yau maximum principle and Yau's Liouville theorem exclude both negative and zero Einstein curvature. We also prove a non-Einstein obstruction under uniformly negative Ricci curvature and obtain compact equivariant continuations.
\end{abstract}
\maketitle
\tableofcontents

\section{Introduction}
A relevant question in the geometry of homolomorphic and isometric immersions asks which K\"ahler--Einstein metrics can be locally induced from the Fubini--Study metric of a complex projective space. In particular, must the Einstein constant of such a metric be positive? Although the Fubini–Study metric has positive curvature, the Gauss equation does not, by itself, impose a sign on the Ricci curvature of a Kähler submanifold.

Calabi's theory \cite{Calabi1953} provides local criteria and a strong uniqueness theorem for K\"ahler immersions (see also the monograph \cite{LoiZedda2018} for a systematic treatment). Hulin subsequently showed that a projectively immersed Einstein germ admits a complete immersed continuation in the same projective space \cite[Proposition~4.5]{Hulin1996} and established positivity for compact projectively immersed K\"ahler--Einstein manifolds \cite{Hulin2000}. More recently, Arezzo, Li and Loi formulated the positivity question explicitly for arbitrary finite codimension \cite[Problem~1.4]{ArezzoLiLoi2025}. Their Corollary~1.7 rules out complete noncompact Ricci-flat K\"ahler submanifolds in finite projective dimension; together with Hulin's continuation and compact positivity theorem, this excludes locally projectively immersed Ricci-flat germs. The case of negative Einstein curvature remains the relevant one in the general question. Our proof under a local symmetry hypothesis will treat both nonpositive Einstein constants by a bounded Hamiltonian potential, independently of the aforementioned Ricci-flat nonimmersion result.

Rigidity results under additional geometric assumptions provide further context. For instance, radial and Ricci-flat settings have been investigated in \cite{ArezzoLoiPlaciniZedda2026,LoiSalisZuddas2018,LoiSalisZuddas2021}, while  K\"ahler--Einstein immersions have been studied in the toric settings in \cite{DiScalaSombra2025,Huang2026,MannoSalis2022,MannoSalis2026,Salis2017}. A recent preprint of Loi and Caib\u{a}r \cite{LoiCaibar2026} classifies compact projectively induced K\"ahler--Einstein surfaces without imposing any additional assumption. These works illustrate the activity around the classification problem. Here we address a different, dimension-independent question: how much infinitesimal symmetry is enough to force the sign of the Einstein constant?

Throughout the paper, a \emph{local holomorphic symmetry} means a real vector field that is simultaneously Killing and real-holomorphic. Our main theorem requires only one nontrivial germ of such a field.
\begin{theorem}[Positivity from one local symmetry]\label{thm:main}
Let $(M,g,J)$ be a connected K\"ahler--Einstein manifold, and suppose that $g$ admits a local K\"ahler immersion into $\CP^N$ for some $N<\infty$. If a nonzero real-holomorphic Killing field is defined on a nonempty open subset of $M$, then the Einstein constant in
$$
 \Ric(g)=\lambda g, \qquad \lambda\in\RR,
$$
satisfies $\lambda>0$.
\end{theorem}

Since K\"ahler--Einstein metrics are real analytic, Calabi's global character of resolvability implies that the existence of a local projective immersion near one point ensures local projective immersibility near every point of the connected manifold (Lemma~\ref{lem:propagation}). Thus, even when their original domains are disjoint, the immersion and the Killing field can be realized on a common neighborhood. The symmetry hypothesis is infinitesimal: neither a global group action nor any prescribed dimension of an isometry group is required. A local Killing field need not have a complete flow on its domain and the original metric need not be complete. Consequently, the local symmetry cannot be treated as arising from a global one-parameter group of isometries, and the classical Bochner vanishing argument for compact manifolds is not directly applicable.
The proof instead makes the projective immersion itself a vehicle for continuing the symmetry.

The first ingredient is a fact about \emph{any} connected K\"ahler manifold with a global full K\"ahler immersion
$$
 F:(M,g,J)\longrightarrow(\CP^N,g_{\FS}).
$$
By Calabi rigidity, a germ of a real-holomorphic Killing field is induced by a unique element of $\mathfrak{pu}(N+1)$. The additional observation is that its ambient projective Killing field, tangent to $F(M)$ on a nonempty open subset of the source, remains tangent along \emph{every} immersed branch, even when $F$ is not injective. Thus every local real-holomorphic Killing field extends uniquely to a global one on $M$ whenever a global full projective immersion is given.

Combining this observation with the complete-continuation theorem of Hulin yields the following global consequence, proved in Section~\ref{sec:Hulin}.
\begin{theorem}[Equivariant continuation of local symmetries]\label{thm:continuation}
Let $(U,g,J)$ be a connected K\"ahler--Einstein germ endowed with a full K\"ahler immersion into $\CP^N$, $N<\infty$. Then there exists a connected complete K\"ahler--Einstein continuation of this immersed germ,
\[
 \widehat F:(\widehat M,\widehat g,\widehat J)\longrightarrow \CP^N,
\]
with the same Einstein constant and a global K\"ahler immersion. Every germ of a real-holomorphic Killing field on $U$ extends uniquely to a global real-holomorphic Killing field on $\widehat M$. The extension is induced by a projective unitary Killing field on $\CP^N$, is complete and has uniformly bounded norm.

The extension identifies the Lie algebra of local real-holomorphic Killing fields with a Lie subalgebra of $\mathfrak{pu}(N+1)$.
\end{theorem}
A noteworthy feature is that the extension takes place directly on $\widehat M$, without passing to a simply connected covering space or imposing any condition on its fundamental group.

The projective realization supplies more than a bounded field: the ambient moment map gives a bounded global Hamiltonian potential $u$ for the extended field. On an Einstein continuation this potential satisfies $\Delta u+2\lambda u=c$. The Omori--Yau maximum principle excludes $\lambda<0$; for $\lambda=0$, it makes $u$ harmonic, and Yau's Liouville theorem forces $u$ to be constant. This proves Theorem~\ref{thm:main}.

The maximum-principle mechanism is not restricted to Einstein metrics. We establish the following complementary result, which separates the intrinsic and extrinsic hypotheses used in the proof.
\begin{theorem}[Negative Ricci obstruction]\label{thm:negative}
Let $(M,g,J)$ be a connected complete K\"ahler manifold admitting a K\"ahler immersion into $\CP^N$, $N<\infty$. Assume that there are constants $C\ge c>0$ such that
$$
 -C g\le \Ric(g)\le -c g.
$$
Then $M$ has no nonzero local real-holomorphic Killing fields. 
\end{theorem}

In fact, the lower Ricci bound is only needed to ensure the weak maximum principle at infinity; we give that more flexible formulation in Theorem~\ref{thm:RicGeneral}. 

Once the sign is known to be positive, the complete Einstein continuation becomes compact by Bonnet--Myers. We record two consequences useful beyond the proof of the sign theorem: the entire local Killing algebra acts by global holomorphic isometries on a compact simply connected projective continuation; and a prescribed compact semisimple symmetry action continues equivariantly, with the dimension of its principal orbits preserved and with an ample linearized polarization (Corollary~\ref{cor:semisimple}). These statements require no assumption on the dimension of the generic orbits. 

It is important to distinguish the ingredients supplied by earlier work from the conclusion established here. 
The local projective-unitary realization follows from Calabi rigidity and the equivariance of the global holomorphic isometry group of an already immersed K\"ahler manifold is recorded, for example, by Zeghib \cite[Theorem~8.1]{Zeghib2016}.
The complete Einstein continuation is Hulin's theorem. The point used here is that a \emph{single local} infinitesimal symmetry descends to a globally defined, bounded field on that continuation. The resulting positivity criterion gives a partial answer to \cite[Problem~1.4]{ArezzoLiLoi2025} under an infinitesimal local symmetry hypothesis.

\medskip
The paper is organized as follows. Section~\ref{sec:local} uses Calabi rigidity to realize germs of real-holomorphic Killing fields by infinitesimal projective unitary transformations. Section~\ref{sec:global} establishes analytic tangency and proves global extension of local symmetries whenever a global full projective immersion is available. Section~\ref{sec:Hulin} applies this extension principle to complete continuations of K\"ahler--Einstein germs and integrates the resulting infinitesimal actions. Section~\ref{sec:Ricci} combines the boundedness of projectively induced Killing fields with the weak maximum principle to obtain an obstruction under negative Ricci curvature. Section~\ref{sec:positivity} proves the positivity criterion and derives compact, homogeneous and equivariantly polarized continuations.

\section{Projective realization of local holomorphic symmetries}\label{sec:local}

We first explain how Calabi rigidity associates to each germ of a real-holomorphic Killing field a unique infinitesimal projective unitary transformation, provided that a full K\"ahler immersion into $\CP^N$ is given.  At this step, neither the Einstein condition nor completeness is involved.

\smallskip
Throughout the paper, we normalize the Fubini--Study form on an affine chart of $\CP^N$ with coordinates $(\zeta_1,\ldots,\zeta_N)$ by
\[
 \omega_{\FS}=\sqrt{-1}\,\partial\bar\partial\log\left(1+\sum_{j=1}^{N}|\zeta_j|^2\right).
\]
A holomorphic isometric immersion will be called a \emph{K\"ahler immersion}. Its image is \emph{full} if it is not contained in a proper projective linear subspace. Every immersion may be made full by replacing the target with its projective linear span. Moreover, by the identity theorem for holomorphic functions, the restriction of a full holomorphic immersion from a connected manifold to any nonempty open subset remains full.

\smallskip
In the Einstein setting, local projective inducedness can be assumed at a single point rather than at every point. Indeed, by the global character of Calabi's resolvability criterion \cite{Calabi1953}, local projective immersibility on a connected real-analytic K\"ahler manifold propagates from one point to every other point, with the minimal full ambient dimension unchanged; see also \cite[Section~1]{HuangLi2025}. We record this classical consequence in the form needed below.

\begin{lemma}\label{lem:propagation}
Let $(M,g,J)$ be a connected real-analytic K\"ahler manifold. If some nonempty open subset admits a K\"ahler immersion into $\CP^N$, then every point of $M$ has a neighborhood admitting such an immersion into the same projective space.
The minimal full ambient dimension is independent of the point.
\end{lemma}

The form of Calabi rigidity needed here is the following \cite{Calabi1953}.
\begin{theorem}[Calabi rigidity]\label{thm:Calabi}
Suppose that $F_1,F_2:M\to\CP^N$ are full K\"ahler immersions of the same connected K\"ahler manifold $M$. There is a unique $A\in\PU(N+1)$ such that $F_2=A\circ F_1$.
\end{theorem}
\begin{remark}\label{rem:full}
A consequence of the uniqueness in Theorem~\ref{thm:Calabi}, used repeatedly below, is that if $F:M\to\mathbb{CP}^N$ is a full K\"ahler immersion of a connected manifold and $B\in\PU(N+1)$ satisfies $B\circ F=F$ on a nonempty open subset of $M$, then $B$ is the identity.
\end{remark}

For a point $p$ of a K\"ahler manifold, let $\Kill^{\mathrm{hol}}_p$ denote the Lie algebra of germs at $p$ of real vector fields $Y$ satisfying
\[
 \mathcal L_Yg=0,\qquad \mathcal L_YJ=0.
\]
Germs of Killing fields form a finite-dimensional Lie algebra. Every real-holomorphic vector field is automatically real analytic, since its $(1,0)$-part has holomorphic coefficients in local holomorphic coordinates.
The metric itself is also real analytic wherever a local K\"ahler immersion exists: indeed, $g=F^*g_{\FS}$, where $F$ is holomorphic and the Fubini--Study metric is real analytic. Independently, K\"ahler--Einstein metrics are real analytic in holomorphic coordinates, by analytic elliptic regularity for the Einstein equation.

\smallskip
For $A\in\mathfrak{pu}(N+1)$ write
\begin{equation}\label{eq:fundamental}
 A^\#_z=\left.\frac{d}{dt}\right|_{t=0}\exp(-tA)\cdot z,
 \qquad z\in\CP^N.
\end{equation}
Throughout the paper, fundamental vector fields of left actions are defined by
\[
\xi^\#_x=
\left.\frac{d}{dt}\right|_{t=0}
\exp(-t\xi)\cdot x.
\]
With this convention, $\xi\mapsto\xi^\#$ is a Lie algebra homomorphism.
\begin{proposition}\label{prop:local}
Let $F:(U,g,J)\to\CP^N$ be a full K\"ahler
immersion, with $U$ connected. For every $p\in U$
there exists a unique injective Lie algebra
homomorphism
$$
\rho_p:\Kill^{\mathrm{hol}}_p
\hookrightarrow\mathfrak{pu}(N+1)
$$
characterized as follows: for each germ $[Y]_p$, the element $A=\rho_p([Y]_p)$ is the unique infinitesimal projective unitary generator whose fundamental vector field is $F$-related to $Y$, that is,
\begin{equation}\label{eq:Frelated}
dF(Y)=A^\#\circ F
\qquad
\text{on a sufficiently small neighborhood of }p.\end{equation}
Moreover, if $Y$ represents a germ at $p$, then $\rho_q([Y]_q)=\rho_p([Y]_p)$ for all
$q$ sufficiently close to $p$.
\end{proposition}
\begin{proof}
Fix a representative $Y$ of the germ. Choose a connected neighborhood $V$ of $p$ on which $Y$ and its local flow $\varphi_t$ are defined, shrinking it so that $\varphi_t(V)\subset U$ for small $t$. Since $\varphi_t$ preserves both $g$ and $J$, the maps $F|_V$ and $F\circ\varphi_t|_V$ are full K\"ahler immersions of the same metric. Theorem~\ref{thm:Calabi} gives a unique $B_t\in\PU(N+1)$ with
\begin{equation}\label{eq:localflow}
 F(\varphi_t(x))=B_t\cdot F(x), \qquad x\in V.
\end{equation}

Choose a sufficiently small connected neighborhood $W$ of $p$ such that $\varphi_s(W)\subset V$ for all sufficiently small $s$. Then, for sufficiently small $t,s$, the local flow identity and \eqref{eq:localflow} give
\[
B_{t+s}F(x)
=F(\varphi_t(\varphi_s(x)))
=B_tB_sF(x),
\qquad x\in W.
\]
Since $F|_W$ is full, Calabi uniqueness yields
$B_{t+s}=B_tB_s$.

We verify the required regularity rather than assuming it.
If $t_j\to0$, compactness of $\PU(N+1)$ ensures that every subsequence of $B_{t_j}$ admits a convergent  subsequence. By \eqref{eq:localflow} and the continuity of the local flow, any such limit fixes $F(V)$ pointwise and hence equals the identity (cf. Remark~\ref{rem:full}). It follows that $B_t\to\id$ as $t\to0$. The local homomorphism property $B_{t+s}=B_tB_s$ then gives
\[
\lim_{s\to0}B_{t+s}
=B_t\,\lim_{s\to0}B_s=B_t
\qquad
\text{for every sufficiently small }t.\]
Thus $t\mapsto B_t$ is continuous in a neighborhood of $0$.
By the smoothness of continuous local homomorphisms between Lie groups, the map $t\mapsto B_t$ is smooth. Recall that every smooth local one-parameter subgroup of a Lie group is locally an exponential curve, uniquely determined by its derivative at the identity. Thus, setting
\[
A=-\left.\frac{d}{dt}\right|_{t=0}B_t
\in\mathfrak{pu}(N+1),
\]
we obtain $B_t=\exp(-tA)$ for all sufficiently small $t$.
Differentiating \eqref{eq:localflow} gives \eqref{eq:Frelated}.

The generator $A$ depends only on the germ $[Y]_p$. Indeed, suppose that $A^\#\circ F=B^\#\circ F$ on some connected neighborhood $W$ of $p$. Then the ambient projective Killing field $(A-B)^\#$ vanishes along $F(W)$, so its flow fixes $F(W)$ pointwise. Since $F|_W$ is full, Calabi rigidity forces this flow to be trivial, and hence $A=B$.

Moreover, the resulting map $\rho_p$ is injective: if $\rho_p([Y]_p)=0$, then \eqref{eq:Frelated} gives $dF(Y)=0$ near $p$. Since $F$ is an immersion, this implies $[Y]_p=0$.

The uniqueness of the ambient generator also implies linearity: applying it to the $F$-related identities for $Y$, $Z$, and $aY+bZ$ gives $\rho_p(aY+bZ)=a\rho_p(Y)+b\rho_p(Z)$. Similarly, $F$-related vector fields have $F$-related Lie brackets, and the sign convention \eqref{eq:fundamental} gives $\rho_p([Y,Z])=[\rho_p(Y),\rho_p(Z)]$.
\end{proof}

Since $\mathfrak{pu}(N+1)$ is a compact Lie algebra, Proposition~\ref{prop:local} has the following immediate consequence.
\begin{corollary}
If a K\"ahler germ admits a full K\"ahler immersion into $\CP^N$, then its local real-holomorphic Killing algebra is an algebra of compact type. In particular, it admits an $\operatorname{ad}$-invariant positive-definite inner product, has dimension at most $(N+1)^2-1$, and contains no noncompact simple real Lie algebra as a subalgebra.
\end{corollary}

\section{Globalizing local symmetries of a projective immersion}\label{sec:global}
The preceding representation associates an ambient projective Killing field with each local symmetry. Its mere existence does not yet guarantee that it restricts to a vector field along the entire immersed manifold. We now prove precisely this tangency statement under the assumption of a \emph{global immersion}, but without completeness or simple connectivity. In general, local Killing fields on real-analytic Riemannian manifolds can develop monodromy (cf. Example~\ref{ex:monodromy}). Nomizu's classical extension theorem \cite{Nomizu1960} assumes simple connectivity to exclude this obstruction. Here a fixed ambient projective generator eliminates monodromy.

\begin{lemma}\label{lem:tangent}
Let $F:M\to\CP^N$ be a holomorphic immersion of a connected complex manifold, and let $Z$ be a real-holomorphic  vector field on $\CP^N$. If $Z$ is tangent to the immersed branch $F(M)$ along $F(V)$ for some nonempty open subset $V\subset M$, then $Z\circ F$ is tangent to the immersed tangent space $dF(TM)$ at every point of $M$.
\end{lemma}
\begin{proof}
Since $F$ is a holomorphic immersion, the quotient
\[
\mathcal N_F=F^*T^{1,0}\CP^N/dF(T^{1,0}M)
\]
is a holomorphic vector bundle on $M$. The $(1,0)$-part of $Z$ induces a holomorphic section of $\mathcal N_F$, whose vanishing is equivalent to tangency along the immersed branch. This section vanishes on $V$, hence everywhere on connected $M$ by the identity theorem. The argument takes place on the source and therefore applies separately to every immersed branch.
\end{proof}

\begin{proposition}[Local-to-global rigidity]\label{prop:global}
Let $F:(M,g,J)\to\CP^N$ be a full K\"ahler immersion of a connected K\"ahler manifold. Every germ of a real-holomorphic Killing field at any point of $M$ admits a unique extension to a global real-holomorphic Killing field on $M$. It is $F$-related to a unique ambient projective unitary Killing field.
\end{proposition}
\begin{proof}
Let $[Y]_p$ be a germ. Proposition~\ref{prop:local} provides a unique $A\in\mathfrak{pu}(N+1)$ such that $dF(Y)=A^\#\circ F$ near $p$. The issue is to show that the same $A$, rather than a succession of locally chosen generators, determines a field on all of $M$. The ambient field $A^\#$ is real analytic. Lemma~\ref{lem:tangent} shows that it is tangent to the immersed branch everywhere on $M$. Since $dF$ is injective, this defines a unique smooth field $\widetilde Y$ on $M$ by
\begin{equation}\label{eq:lift}
 dF(\widetilde Y)=A^\#\circ F.
\end{equation}
The resulting field is smooth because $dF(TM)$ is a smooth subbundle and $dF$ is a bundle isomorphism onto its image. On each sufficiently small domain where $F$ is an embedding, $\widetilde Y$ is the restriction of the ambient field $A^\#$ to the embedded complex submanifold. It is therefore real-holomorphic and Killing: its local flow is intertwined by $F$ with the holomorphic isometric flow of $A^\#$. These local conclusions agree on overlaps and hence hold globally, even when distinct points of $M$ have the same image under $F$. Equation \eqref{eq:lift} agrees with $Y$ near $p$. Uniqueness follows because a global Killing field vanishing on a nonempty open subset of a connected Riemannian manifold vanishes identically.
\end{proof}
\begin{example}\label{ex:monodromy}
The existence of a global projective immersion in Proposition~\ref{prop:global} cannot be replaced by local projective inducedness. Let $n\ge2$ and equip $\Omega=\mathbb C^n\setminus\{0\}$ with the restriction of the Fubini--Study metric under the standard affine embedding into $\CP^n$. The involution $\sigma(z)=-z$ is a fixed-point-free holomorphic isometry, so the quotient $M=\Omega/\langle\sigma\rangle$ carries a K\"ahler--Einstein metric which is locally projectively induced.

Let $S=\operatorname{diag}(1,-1,\ldots,-1)$ represent $\sigma$ on $\CP^n$, and choose $A=E_{01}-E_{10}\in\mathfrak{su}(n+1)\cong\mathfrak{pu}(n+1)$, where $E_{ij}$ denotes the standard matrix unit. 
Since $SAS^{-1}=-A$, the Killing field $Y=A^\#|_\Omega$ satisfies $\sigma_*Y=-Y$.
It induces a nonzero local real-holomorphic Killing field on $M$, but this field cannot extend globally. Indeed, such an extension would lift to a $\sigma$-invariant Killing field on $\Omega$, which would coincide with $Y$ by unique continuation, a contradiction.
\end{example}
The uniqueness assertion in Proposition~\ref{prop:global} immediately yields the following description of the local Killing algebras.
\begin{corollary}\label{cor:sheaf}
Under the hypotheses of Proposition~\ref{prop:global}, restriction to germs induces Lie algebra isomorphisms
\[
 \Kill^{\mathrm{hol}}(M)\xrightarrow{\ \simeq\ }\Kill^{\mathrm{hol}}_p
 \qquad(p\in M),
\]
where $\Kill^{\mathrm{hol}}(M)$ denotes the global real-holomorphic Killing fields. In particular, every local real-holomorphic Killing field extends around any closed curve without monodromy, and the dimensions of the germ algebras are constant on $M$.
\end{corollary}

\begin{proposition}\label{prop:bounded}
In the setting of Proposition~\ref{prop:global}, every global field $\widetilde Y$ extending a local real-holomorphic Killing field satisfies
\begin{equation}\label{eq:uniformbound}
 \sup_M |\widetilde Y|_g
 \le\sup_{z\in\CP^N}|A^\#_z|_{g_{\FS}}<\infty.
\end{equation}
If $(M,g)$ is complete, then $\widetilde Y$ is a complete vector field. Its flow $\varphi_t$ is $F$-equivariant:
\begin{equation}\label{eq:equivflow}
 F\circ\varphi_t=\exp(-tA)\circ F,\qquad t\in\RR.
\end{equation}
\end{proposition}
\begin{proof}
Since $F$ is isometric, \eqref{eq:lift} gives the pointwise identity
$ |\widetilde Y_x|_g=|A^\#_{F(x)}|_{g_{\FS}}$.
The continuous ambient norm is bounded on compact $\CP^N$, proving \eqref{eq:uniformbound}. A bounded smooth vector field on a complete Riemannian manifold is complete: any integral curve defined on a finite time interval stays in a bounded metric ball, which is compact by Hopf--Rinow, and can therefore be extended. Finally, for each $p\in M$, both curves $t\mapsto F(\varphi_t(p))$ and $t\mapsto \exp(-tA)\cdot F(p)$ are integral curves of $A^\#$ starting at $F(p)$. The equivariance identity follows by uniqueness.
\end{proof}

\section{Equivariant complete continuation in the Einstein case}\label{sec:Hulin}
The globalization result of Section~\ref{sec:global} assumes that the immersion is defined on the entire source. An Einstein germ need not come with such an immersion globally; Hulin's continuation theorem supplies exactly the missing geometric input. Hulin's complete-continuation theorem \cite[Proposition~4.5]{Hulin1996} provides, starting from a sufficiently small immersed Einstein germ, a connected complete K\"ahler--Einstein manifold $\widehat M$ with a K\"ahler immersion $\widehat F$ into the same finite-dimensional projective space, extending the germ.
This continuation theorem is stated in complex dimension at least two. In complex dimension one, the source has constant Gaussian curvature. The local classification of K\"ahler immersions between complex space forms by Nakagawa and Ogiue \cite[Theorem~1]{NakagawaOgiue1976}, together with Calabi rigidity \cite{Calabi1953}, identifies its full projective immersion locally with a Veronese immersion of $\CP^1$. This gives the required complete model in dimension one.

\begin{lemma}\label{lem:curve}
Let $(\Sigma,g)$ be a connected K\"ahler--Einstein complex
curve. If some nonempty open subset admits a K\"ahler
immersion into $\CP^N$, $N<\infty$, then its Gaussian
curvature is positive. Moreover, every sufficiently small
immersed germ extends to a complete projectively immersed
model isometric to $\CP^1$ endowed with a positive integer
multiple of the Fubini--Study metric.
\end{lemma}

We may therefore use the continuation statement below in every complex dimension, treating dimension one by Lemma~\ref{lem:curve}.

\begin{theorem}\label{thm:equiv}
Let $(U,g,J)$ be a K\"ahler--Einstein germ with a full K\"ahler immersion $F:U\to\CP^N$, $N<\infty$. Let
\[
 \widehat F:(\widehat M,\widehat g,\widehat J)\longrightarrow\CP^N
\]
be a connected complete continuation of this immersed germ. Then the algebra $\Kill^{\mathrm{hol}}_p$, for $p$ in the original germ, is canonically isomorphic to the algebra of global real-holomorphic Killing fields of $\widehat M$. Every such global field is complete and has bounded norm. There is an injective Lie algebra homomorphism
\[
 \widehat\rho:\Kill^{\mathrm{hol}}_p\longrightarrow\mathfrak{pu}(N+1)
\]
which intertwines the infinitesimal action on $\widehat M$ with the projective unitary action on the ambient space.
\end{theorem}
\begin{proof}
The restriction of $\widehat F$ to the original germ is full; therefore $\widehat F$ is itself full. Proposition~\ref{prop:global} applied to the connected globally immersed manifold $\widehat M$ extends every germ. The identification of algebras follows from Corollary~\ref{cor:sheaf}. The representation is the map furnished by Proposition~\ref{prop:local}, with the same projective generator by uniqueness. Completeness, boundedness and equivariance follow from Proposition~\ref{prop:bounded}.
\end{proof}
Together with the existence of the complete continuation recalled above (and Lemma~\ref{lem:curve} in complex dimension one), this proves Theorem~\ref{thm:continuation} stated in the Introduction.

\begin{corollary}
Suppose a finite-dimensional real Lie algebra $\mathfrak h$ acts effectively on the immersed Einstein germ by real-holomorphic Killing fields. Then the infinitesimal action extends to $\widehat M$ and integrates to a global action of the simply connected Lie group $H_{\mathrm{sc}}$ with Lie algebra $\mathfrak h$. Moreover, there is a smooth homomorphism
\[
 \rho:H_{\mathrm{sc}}\longrightarrow\PU(N+1)
\]
for which $\widehat F$ is equivariant. The action has discrete kernel.
\end{corollary}
\begin{proof}
Theorem~\ref{thm:equiv} supplies an injective Lie algebra homomorphism from $\mathfrak h$ into the complete global Killing fields of $\widehat M$ and a compatible homomorphism into $\mathfrak{pu}(N+1)$. By the Lie--Palais integration theorem \cite{Palais1957}, the former integrates to an action of a simply-connected group $H_{\mathrm{sc}}$. The latter integrates by simple connectivity to a homomorphism into $\PU(N+1)$. Infinitesimal equivariance implies equivariance for the connected group. An infinitesimally faithful action has discrete kernel.
\end{proof}

The argument works on the continuation $\widehat M$ itself, without passing to its universal cover: the fixed projective generator ensures that the extended fields descend through any monodromy. In particular, neither compactness nor simple connectivity is used before the sign of the Einstein constant is established.

\section{Negative Ricci curvature and the maximum principle}\label{sec:Ricci}
The projective argument has now supplied a bounded \emph{global} Killing field. We turn to the analytic obstruction to its existence in negative Ricci curvature. The argument needs only a weak maximum principle for bounded-above functions; in particular, no Einstein condition is required in this section.

We use the sign convention $\Delta=\operatorname{div}\nabla$ for the Laplace--Beltrami operator. Recall that the \emph{weak maximum principle at infinity} holds on a complete Riemannian manifold if, for every $u\in C^2(M)$ bounded above with $S=\sup_M u$, there is a sequence $x_j\in M$ such that
\begin{equation}\label{eq:WMP}
 u(x_j)>S-\frac1j,\qquad \Delta u(x_j)<\frac1j.
\end{equation}
A complete Riemannian manifold with Ricci curvature bounded from below has this property by the Omori--Yau maximum principle \cite{Yau1975}. The weaker formulation \eqref{eq:WMP} is all we need.

\begin{lemma}\label{lem:Bochner}
Let $(M,g)$ be a connected complete Riemannian manifold satisfying the weak maximum principle at infinity. If
\[
 \Ric(g)\le -c g\qquad(c>0),
\]
then every bounded global Killing field on $M$ is identically zero.
\end{lemma}
\begin{proof}
Let $Y$ be such a field. For a Killing field, the Bochner identity, with our sign convention, reads
$$
 \frac12\Delta |Y|^2=|\nabla Y|^2-\Ric(Y,Y).
$$
Consequently,
\begin{equation}\label{eq:deltau}
 \Delta |Y|^2\ge 2c\,|Y|^2.
\end{equation}
Since $Y$ is bounded, $S=\sup_M |Y|^2<\infty$. If $S>0$, the weak maximum principle at infinity provides points  $x_j$ satisfying \eqref{eq:WMP} with $u= |Y|^2$. For large $j$, $|Y|^2(x_j)>S/2$ and hence
\[
 \frac1j>\Delta |Y|^2(x_j)\ge 2c\,|Y|^2(x_j)>cS,
\]
a contradiction. Thus $S=0$ and $Y=0$.
\end{proof}

\begin{theorem}\label{thm:RicGeneral}
Let $(M,g,J)$ be a connected complete K\"ahler manifold with a K\"ahler immersion into $\CP^N$, $N<\infty$. Suppose the weak maximum principle at infinity holds and
\[
 \Ric(g)\le -c g
\]
for some $c>0$. Then every local real-holomorphic Killing field on $M$ vanishes. In particular, the conclusion holds under the two-sided bound
\[
 -C g\le\Ric(g)\le-c g,\qquad C\ge c>0.
\]
\end{theorem}
\begin{proof}
Pass to the minimal projective linear subspace containing the immersed image so that the global immersion $F$ is full. By Proposition~\ref{prop:global}, any germ extends to a global real-holomorphic Killing field $\widetilde Y$. Proposition~\ref{prop:bounded} gives $\sup_M|\widetilde Y|<\infty$. The field vanishes by Lemma~\ref{lem:Bochner}, and so does the initial germ. The two-sided Ricci estimate implies the weak maximum principle by the Ricci lower bound.
\end{proof}
In particular, the two-sided Ricci bound gives Theorem~\ref{thm:negative} stated in the Introduction.

\begin{remark}
This theorem makes no Einstein assumption and the second fundamental form of the immersion is unrestricted. The uniform negative upper bound is used only in \eqref{eq:deltau}; the weak maximum principle controls bounded functions at infinity. Completeness and the existence of a \emph{global} projective immersion are required here. In contrast, the Einstein case in the next section will require only a local projective immersion, because Hulin's theorem supplies a global one on a complete continuation.
\end{remark}

\section{Positivity and compact equivariant continuations}\label{sec:positivity}
We now use the ambient projective moment map to obtain a bounded global Hamiltonian potential for each Killing field on a complete projectively immersed K\"ahler manifold. In the Einstein case, its Laplace equation treats negative and zero Einstein constants by the same mechanism.

\begin{proposition}[Bounded Hamiltonian obstruction]\label{prop:boundedmoment}
Let $(M,g,J)$ be a connected complete K\"ahler--Einstein manifold with $\Ric(g)=\lambda g$ and $\lambda\le0$. If a real-holomorphic Killing field $Y$ has a globally defined bounded Hamiltonian potential $u$, then $Y=0$.
\end{proposition}
\begin{proof}
Let $\omega(\cdot,\cdot)=g(J\cdot,\cdot)$ be the K\"ahler form. We choose the sign of $u$ so that
\[
du=\iota_Y\omega,\qquad Y=-J\nabla u.
\]
Recall that a Killing field satisfies $\nabla^*\nabla Y=\Ric^\sharp(Y)$, whereas
\[
 \nabla^*\nabla(\nabla u)=-\nabla(\Delta u)-\Ric^\sharp(\nabla u).
\]
Since $J$ is parallel and the Ricci endomorphism commutes with $J$, these identities yield
\[
 \nabla(\Delta u)=-2\Ric^\sharp(\nabla u)=-2\lambda\nabla u.
\]
Thus, for some constant $c\in\RR$,
\begin{equation}\label{eq:momentlaplace}
 \Delta u+2\lambda u=c.
\end{equation}

If $\lambda<0$, the bounded function $v=u-c/(2\lambda)$ satisfies $\Delta v=-2\lambda v$. The Omori--Yau weak maximum principle, valid because $M$ is complete with Ricci curvature bounded below, applied to $v$ and $-v$ gives $\sup_M v\le0$ and $\inf_M v\ge0$. Hence $v=0$.

If $\lambda=0$, equation \eqref{eq:momentlaplace} gives $\Delta u=c$. Applying the same maximum principle to $u$ and $-u$ shows that $c=0$. Thus $u$ is bounded and harmonic. Since $\Ric(g)=0$, Yau's Liouville theorem for bounded harmonic functions on complete manifolds with nonnegative Ricci curvature \cite[Corollary~1]{Yau1975} implies that $u$ is constant. In either case $Y=-J\nabla u=0$.
\end{proof}

The projective setting supplies precisely the bounded potential required in Proposition~\ref{prop:boundedmoment}. If $F:(M,g,J)\to\CP^N$ is a full K\"ahler immersion and $Y$ is a global real-holomorphic Killing field, let $A\in\mathfrak{su}(N+1)$ represent its unique ambient infinitesimal generator, with the convention that its flow is $\exp(-tA)$. Put $H=-\sqrt{-1}A$, which is Hermitian. The standard projective moment map component
\[
 \mu_A([Z])=\frac{Z^*HZ}{Z^*Z},\qquad Z\in\CC^{N+1}\setminus\{0\},
\]
satisfies $d\mu_A=\iota_{A^\#}\omega_{\FS}$. Consequently,
\begin{equation}\label{eq:globalmoment}
 u=\mu_A\circ F,\qquad du=\iota_Y\omega,
\end{equation}
and $u$ is globally bounded, since its values lie between the smallest and largest eigenvalues of $H$.

\begin{proof}[Proof of Theorem~\ref{thm:main}]
Let $Y$ be a nonzero local real-holomorphic Killing field. By Lemma~\ref{lem:propagation}, we may restrict to a connected germ on which both $Y$ and a full K\"ahler immersion into $\CP^N$ are defined. Hulin's theorem (or Lemma~\ref{lem:curve} in complex dimension one) furnishes a connected complete K\"ahler--Einstein continuation with the same Einstein constant $\lambda$ and a global projective immersion.

By Theorem~\ref{thm:equiv}, $Y$ extends to a nonzero global real-holomorphic Killing field on this continuation. Equation~\eqref{eq:globalmoment} provides a bounded global Hamiltonian potential. Proposition~\ref{prop:boundedmoment} excludes $\lambda\le0$, and hence $\lambda>0$.
\end{proof}
\begin{remark}
The negative Einstein case also excludes nonzero
ambient holomorphic projective vector fields tangent
to a full immersed germ, without assuming that they
are Killing. More precisely, let
$F:(U,g)\to\mathbb{CP}^N$ be a full K\"ahler
immersion, with $U$ connected and
$\Ric(g)=\lambda g$, $\lambda<0$.
Let $A\in\mathfrak{pgl}(N+1,\mathbb C)$, and suppose
that its fundamental field $A^\#$ satisfies
\[
A^\#(F(p))\in dF_p(T_pU)
\qquad\text{for all }p\in U.
\]
The case of complex dimension one is excluded by
Lemma~\ref{lem:curve}. Otherwise, Hulin's complete
continuation and Lemma~\ref{lem:tangent} yield a global
real-holomorphic vector field $Y$ on a complete Einstein
continuation, induced by the same ambient field $A^\#$.
The field $Y$ is bounded, since the immersion is
isometric and $\mathbb{CP}^N$ is compact.
Its $(1,0)$-part $V=Y^{1,0}$ satisfies the
holomorphic Bochner identity
\[
\frac12\Delta|V|^2
=|\nabla^{1,0}V|^2-\lambda|V|^2
\ge -\lambda|V|^2.
\]
Since $\lambda<0$, the Omori--Yau maximum principle
implies $V=0$, hence $Y=0$. Thus $A^\#$ vanishes on $F(U)$.
Since $F(U)$ is connected, it lies in the
projectivization of a single eigenspace of $A$.
Fullness then forces $A=0$.
\end{remark}

\paragraph{\textbf{Compact and equivariant consequences}}
The positivity criterion is not merely a statement about a curvature constant. Together with Hulin's theorem, it turns the complete continuation into a compact geometric object on which every original local symmetry acts globally. The following statements record the complete package, first for arbitrary local Killing algebras and then for an assigned compact semisimple group.

\begin{corollary}\label{cor:compactcont}
Let $(M,g,J)$ be a connected K\"ahler--Einstein manifold admitting a local K\"ahler immersion into a finite-dimensional projective space, and assume that a nonzero local real-holomorphic Killing field exists. Let $N$ be the smallest integer such that some nonempty open subset of $M$ admits a K\"ahler immersion into
$\mathbb{CP}^N$. Then every sufficiently small immersed germ $(U,g|_U,J|_U)$ admits an open holomorphic isometric embedding into a compact simply connected K\"ahler--Einstein manifold $(\widetilde M,\widetilde g,\widetilde J)$ with a full global K\"ahler immersion $\widetilde F:\widetilde M\to\CP^N$. The Einstein constant agrees with that of $g$. Every local real-holomorphic Killing field on $U$ extends uniquely to a complete global real-holomorphic Killing field on $\widetilde M$, and this extension is induced by a projective unitary Killing field. Moreover, $\widetilde M$ is smooth and projective and $\widetilde F^*\mathcal O_{\CP^N}(1)$ is ample.
\end{corollary}
\begin{proof}
By Theorem~\ref{thm:main}, $\lambda>0$. By Lemma~\ref{lem:propagation}, the minimal full projective dimension is independent of the chosen germ. Choose a connected immersed germ with this minimal full dimension $N$ and apply Hulin's complete-continuation theorem to obtain $(\widehat M,\widehat g,\widehat J,\widehat F)$. Bonnet--Myers implies that $\widehat M$ is compact and its fundamental group finite. Let $\pi:\widetilde M\to\widehat M$ be its universal Riemannian cover, endowed with the pulled-back complex structure and metric and the immersion $\widetilde F=\widehat F\circ\pi$. Then $\widetilde M$ is compact and simply connected, and after shrinking $U$ the original germ embeds holomorphically and isometrically into $\widetilde M$. Its immersion is the restriction of $\widetilde F$ (up to an ambient projective unitary isometry). Theorem~\ref{thm:equiv}, applied to $\widetilde M$, extends all local real-holomorphic Killing fields and realizes them projectively. Finally $\widetilde F^*\mathcal O_{\CP^N}(1)$ is a positive holomorphic line bundle by the isometry condition, hence ample; Kodaira's theorem makes $\widetilde M$ projective.
\end{proof}

\begin{corollary}[Equivariant polarized continuation for semisimple symmetries]\label{cor:semisimple}
Let $(M,g,J)$ be a connected K\"ahler--Einstein manifold admitting a local K\"ahler immersion into a finite-dimensional projective space. Suppose a nontrivial compact connected semisimple Lie group $K$ acts almost effectively on $M$ by holomorphic isometries, and let $\widetilde K\to K$ be the simply connected covering group. Then every sufficiently small immersed germ $(U,g|_U,J|_U)$ has an open holomorphic isometric embedding $\iota:U\to \widetilde M$ into a compact simply connected K\"ahler--Einstein manifold $(\widetilde M,\widetilde g,\widetilde J)$ such that:
\begin{enumerate}[label=\textup{(\roman*)},leftmargin=*]
\item $\widetilde K$ acts almost effectively on $\widetilde M$ by holomorphic isometries; its infinitesimal action on $\iota(U)$ agrees, under $d\iota$, with the infinitesimal $K$-action on $U$;
\item there is a full $\widetilde K$-equivariant K\"ahler immersion $\widetilde F:\widetilde M\to\CP^N$, with $N$ the minimal projective dimension of the original germ and $\widetilde F\circ\iota$ equal to the germ immersion, up to an ambient projective unitary motion;
\item $\widetilde M$ is a smooth projective manifold; writing $G=\widetilde K^{\CC}$ for the simply connected complexification, the holomorphic action extends to $G$, the immersion $\widetilde F$ is $G$-equivariant for a representation $G\to\operatorname{SL}(N+1,\CC)$, and the ample line bundle $L=\widetilde F^*\mathcal O_{\CP^N}(1)$ is $G$-linearized;
\item the maximum dimensions of the $K$-orbits on $M$ and the $\widetilde K$-orbits on $\widetilde M$ are equal. In particular, the principal-orbit codimension is preserved.
\end{enumerate}
\end{corollary}
\begin{proof}
The infinitesimal action of $\mathfrak k=\Lie(K)$ is faithful: a fundamental field vanishing on a nonempty open subset of the connected manifold $M$ vanishes identically, and almost effectiveness then makes its generator zero. Thus Theorem~\ref{thm:main} gives $\lambda>0$. Apply Corollary~\ref{cor:compactcont} to obtain $(\widetilde M,\widetilde g,\widetilde J,\widetilde F,\iota)$, where $\widetilde M$ is compact and simply connected and $\widetilde F$ is full.

Theorem~\ref{thm:equiv} extends the fundamental fields from $U$ to complete real-holomorphic Killing fields on $\widetilde M$, preserving their brackets and the infinitesimal action on $\iota(U)$. Lie--Palais integration gives an action of $\widetilde K$ on $\widetilde M$. The corresponding faithful Lie algebra homomorphism into $\mathfrak{pu}(N+1)$ integrates to a homomorphism $\rho:\widetilde K\to\PU(N+1)$. Infinitesimal equivariance, followed by connectedness of $\widetilde K$, gives $\widetilde F(k\cdot x)=\rho(k)\widetilde F(x)$. The action has discrete kernel by infinitesimal faithfulness, hence finite kernel because $\widetilde K$ is compact. This proves (i) and (ii).

The pullback of the hyperplane bundle has positive Chern curvature $\widetilde\omega=\widetilde F^*\omega_{\FS}$, so it is ample. Since $\widetilde K$ is simply connected, $\rho$ lifts through the finite covering $\operatorname{SU}(N+1)\to\PU(N+1)$ to a unitary representation $\widetilde\rho:\widetilde K\to\operatorname{SU}(N+1)$. The compact complex manifold $\widetilde M$ has complete holomorphic vector fields, so the holomorphic infinitesimal $\mathfrak k$-action complexifies and integrates to a holomorphic action of the simply connected complex group $G=\widetilde K^{\CC}$. Likewise $\widetilde\rho$ complexifies to a holomorphic representation $\widetilde\rho_{\CC}:G\to\operatorname{SL}(N+1,\CC)$. The two holomorphic maps $G\to\CP^N$, $h\mapsto\widetilde F(h\cdot x)$ and $h\mapsto\widetilde\rho_{\CC}(h)\widetilde F(x)$, agree on the real form $\widetilde K$, hence on all of $G$ by the identity principle. Thus $\widetilde F$ is $G$-equivariant, and the pullback of the standard linearization of $\mathcal O_{\CP^N}(1)$ gives a holomorphic $G$-linearization of $L$. This linearization is also algebraic. Since $L$ is ample, some positive tensor power
$L^{q}$ is very ample and hence defines an embedding of $\widetilde M$ into projective space. The induced finite-dimensional holomorphic representation of the complex semisimple group $G$ on $H^0(\widetilde M,L^q)$ is algebraic, so the resulting equivariant embedding makes the $G$-action on $\widetilde M$ algebraic. Since $\widetilde F$ is a holomorphic map between projective varieties, it is algebraic; pulling back the standard algebraic linearization of $\mathcal O_{\CP^N}(1)$ proves (iii).

For (iv), let $r$ be the maximal rank of the fundamental fields on $M$. This rank occurs on a dense open subset, which every nonempty germ meets. The fields on $\widetilde M$ therefore attain rank $r$ at some point of $\iota(U)$. Every $(r+1)$-fold wedge of extended fundamental fields is a real-analytic tensor on $\widetilde M$ vanishing on $\iota(U)$; it must vanish identically because $\widetilde M$ is connected. The maximum rank on $\widetilde M$ is consequently exactly $r$. For compact group actions the maximum orbit dimension equals the dimension of the principal orbits, proving (iv).
\end{proof}

%

A short application illustrates the strength of continuing the \emph{entire} local Killing algebra. If its evaluations span the tangent space at one point, transitivity becomes global on the compact continuation. The resulting homogeneous classification is classical; what is relevant here is to show that \emph{local} transitivity, without completeness of the original germ, propagates to the compact continuation.

\begin{corollary}
Suppose, in addition to the assumptions of Theorem~\ref{thm:main}, that the germs of real-holomorphic Killing fields span $T_pM$ at some point $p\in M$. Then every connected complete projectively immersed continuation of a sufficiently small germ centered at $p$ is a compact homogeneous K\"ahler--Einstein manifold. Its universal cover is a compact simply connected homogeneous K\"ahler manifold and hence a generalized flag manifold, equipped with a homogeneous K\"ahler--Einstein metric.
\end{corollary}
\begin{proof}
By Theorem~\ref{thm:equiv}, the local fields extend to global holomorphic Killing fields on the compact $\widehat M$. Let $K$ be the identity component of its group of holomorphic isometries. This is a compact Lie group. At the chosen point its fundamental fields span the tangent space, so the corresponding $K$-orbit is open. Compactness of $K$ makes all its orbits closed, hence the orbit is simultaneously open and closed. Connectedness of $\widehat M$ gives transitivity. On the compact simply connected universal cover $\widetilde M$, the global Killing fields lift to global real-holomorphic Killing fields spanning the tangent space at a point above the chosen one. The identity component of the holomorphic isometry group of $\widetilde M$ is compact and has an open orbit, hence acts transitively by the same open-and-closed orbit argument. The Borel--Remmert structure theorem \cite{BorelRemmert1962} identifies a compact simply connected homogeneous K\"ahler manifold with a generalized flag manifold. The lifted metric is homogeneous and K\"ahler--Einstein.
\end{proof}

\paragraph{\textbf{Sharpness and boundary cases}}
We conclude with two observations concerning the sharpness of our results.

First, the positivity conclusion in Theorem~\ref{thm:main} is sharp in the sense that no uniform positive lower bound for the Einstein constant can be imposed, even in fixed complex dimension. Indeed, for every positive integer $m$, the metric $g_m=m g_{\mathrm{FS}}$ on $\CP^n$ is K\"ahler--Einstein and projectively induced
by the degree-$m$ Veronese embedding. It admits nonzero real-holomorphic Killing fields, while its Einstein constant is
\[
\lambda_m=\frac{n+1}{m}\longrightarrow0^+.
\]
Here the ambient projective dimension is $\binom{n+m}{m}-1$ and increases with $m$.

Second, finite projective dimension is essential. Indeed, the holomorphic immersion
\[
F:\mathbb D\longrightarrow\CP^\infty,
\qquad
F(z)=[1:z:z^2:\cdots],
\]
induces the K\"ahler form
\[
F^*\omega_{\mathrm{FS}}
=
-\sqrt{-1}\,\partial\bar\partial
\log(1-|z|^2).
\]
The resulting metric is complete and K\"ahler--Einstein, with $\Ric(g)=-2g$, and admits a nonzero real-holomorphic Killing field generated by rotations. Thus the positivity conclusion fails in infinite projective dimension. Moreover, this Killing field has unbounded norm, showing why the boundedness argument used in finite dimension no longer applies.

\bigskip
\paragraph{\textbf{Acknowledgements.}}
The author acknowledges support from the ``Starting Grant'' under contract no.~53\_RSG22SALFIL. He is a member of GNSAGA of the Istituto Nazionale di Alta Matematica (INdAM).

\end{document}